\documentclass[11pt]{amsart}
\title[Cohomology of augmentation ideals in $\ell^1(G)$]{Simplicial coho\-mology of augmentation ideals in~$\ell^1(G)$}

\RequirePackage{amsmath, amssymb, amsfonts, amsthm}

\author{Yemon Choi}

\date{Last updated August 2008.}

\usepackage{concrete} 

\newcounter{pulse}
\numberwithin{pulse}{section}

\RequirePackage{hyperref,verbatim}  

\address{Department of Mathematics, University of Manitoba, Winnipeg, Manitoba R3T 2N2, Canada.}
\address{Current address: D\'epartement des math\'ematiques et de statistique, Universit\'e Laval, Quebec G1V 0A6, Canada.}
\email{y.choi.97@cantab.net}

\theoremstyle{plain}
\newtheorem{thm}[pulse]{Theorem}
\newtheorem*{thm*}{Theorem}

\newtheorem{lem}[pulse]{Lemma}
\newtheorem{coroll}[pulse]{Corollary}
\theoremstyle{definition}
\newtheorem{defn}[pulse]{Definition}
\newtheorem*{notn}{Notation}
\newtheorem{eg}[pulse]{Example}
\theoremstyle{remark}
\newtheorem*{rem}{Remark}
\newtheorem*{unsolved}{Question}

\RequirePackage[PS,nohug]{diagrams}

\newcommand{\xlarr}[1]{{\lTo^{#1}}}
\newcommand{\xrarr}[1]{{\rTo^{#1}}}
\newcommand{\larr}{{\lTo}}
\newcommand{\rarr}{{\rTo}}

\newcommand{\pcat}[1]{{\sf {#1}}} 

\newcommand{\Cst}[2]{
 {\raise0.7ex\hbox{${#1}$} \!\mathord{\left/
 {\vphantom {{#1} {#2}}}\right.\kern-\nulldelimiterspace}
 \!\lower0.7ex\hbox{${#2}$}}
}

\newcommand{\sid}{\mathop{\sf id}\nolimits}
\newcommand{\mc}{\mathcal}
\newcommand{\veps}{\varepsilon}
\newcommand{\blank}{\underline{\quad}}
\newcommand{\iso}{\cong}
\newcommand{\miso}{\underset{1}{\iso}} 
\newcommand{\Ind}{{\mathbb I}} 

\newcommand{\tveps}{{\widetilde{\varepsilon}}}
\newcommand{\td}{{\widetilde{d}}}

\DeclareMathOperator{\Coker}{coker}

\newcommand{\st}{\;:\;}
\newcommand{\defeq}{:=}
\newcommand{\dt}[1]{{\it #1}\/}  
\renewcommand{\emph}[1]{{\sl #1\/}}  

\newcommand{\Ban}{\pcat{Ban}}

\newcommand{\unMod}[2]{{}_{#1}\pcat{unmod}_{#2}} 
\newcommand{\LunMod}[1]{\unMod{#1}{}} 
\newcommand{\lHom}[1]{{}_{#1}\mathop{\rm Hom}\nolimits} 
\newcommand{\cU}{{\mathcal U}}  

\newcommand{\norm}[1]{\Vert{#1}\Vert}
\newcommand{\Lin}[2]{{\mathcal L}({#1},{#2})} 

\newcommand{\tp}{\mathop{\scriptstyle\otimes}}
\newcommand{\ptp}{\mathop{\scriptstyle\widehat{\otimes}}}
\newcommand{\tpR}[1]{\mathop{\underset{#1}{\scriptstyle\otimes}}}
\newcommand{\ptpR}[1]{\mathop{\underset{#1}{\scriptstyle\widehat{\otimes}}}}

\newcommand{\Cplx}{\mathbb C}
\newcommand{\Nat}{\mathbb N}
\newcommand{\Z}{{\mathbb Z}}

\newcommand{\lp}[2][]{\ell^{#2}_{#1}}
\newcommand{\lpsum}[1]{\mathop{\overset{[#1]}{\bigoplus}}}

\newcommand{\id}[1][]{{\sf 1}_{#1}} 
\newcommand{\fu}[1]{{#1}^\#}    

\newcommand{\AI}[1]{I_0(#1)} 

\newcommand{\Cl}[1]{{\mathfrak{Cl}}_{#1}}  

\newcommand{\twist}[1]{{{#1}^\circ}}  

\newcommand{\Ext}{\mathop{\rm Ext}\nolimits}
\newcommand{\Co}[3][]{{\mathcal {#2}^{#3}_{#1}}} 
\newcommand{\dif}{\delta} 

\newcommand{\ct}{com\-muta\-tive-trans\-itive}
\newcommand{\Gronbaek}{Gr{\o}nb{\ae}k}

\begin{document}

\begin{abstract}
Let $G$ be a discrete group. We give a decomposition theorem for the Hochschild cohomology of $\lp{1}(G)$ with coefficients in certain $G$-modules. Using this we show that if $G$ is \ct, the canonical inclusion of bounded coho\-mology of $G$ into simplicial coho\-mology of $\lp{1}(G)$ is an iso\-morphism.

\medskip

\end{abstract}

\maketitle

\begin{section}{Introduction}
The bounded coho\-mology of a (discrete) group $G$ is known to embed as a summand in the simplicial coho\-mology of the convolution algebra $\lp{1}(G)$.
Consequently, knowing that the bounded coho\-mology of $G$ is non-zero, or non-Hausdorff, immediately implies that the simplicial coho\-mology of $\lp{1}(G)$ is non-zero or non-Hausdorff respectively.

In this article we observe that for a wide class of discrete groups, including all torsion-free hyperbolic groups, this summand is the only
non-zero contribution to simplicial coho\-mology: more precisely, the aforementioned inclusion of bounded coho\-mology into simplicial coho\-mology is an isomorphism. The precise statement is given as Theorem~\ref{t:augideal_coho_triv} below.
By standard homological arguments
(see Lemma~\ref{l:simptrivcondn} below)
 we may recast our result as saying that the augmentation ideals for these groups are simplicially trivial, in the sense that the `naive' Hochschild coho\-mology groups $\Co{H}{*}(I_0(G), I_0(G)')$ vanish: see Corollary~\ref{c:final_cor}.
 Thus our work is a partial generalization of results of \Gronbaek\ and Lau (\cite{GronLau_AI}) on weak amenability of such ideals.

Our work is also motivated by the preprint \cite{PW_bddcoho}, in which a version of our decomposition theorem is given for second-degree cohomology: there, the conclusion is stronger because the second bounded coho\-mology of \emph{any} discrete group is known to be a Banach space (no such general result is true for degrees 3 and above).
\begin{rem}
After the main work of this article was completed, the author learned of the article \cite{Pou_gpalg}. The two articles do not overlap much, but the proof of our main decomposition theorem could be significantly shortened if \cite[Corollary 3.7]{Pou_gpalg} were valid. However, the claim made in that corollary seems to require further justification: see Example \ref{eg:SNIPER} below for more details.
\end{rem}
\end{section}

\begin{section}{Notation and homological background}
\label{s:prelim}
Throughout we shall denote the identity map on a Banach space, module or algebra by $\sid$ (it will be clear from context what the domain of $\sid$~is).
Isometric linear isomorphism of two Banach spaces $E$ and $F$ will be denoted by $E\miso F$\/; the dual of a Banach space $E$ will be denoted by~$E'$\/.
Given a family $(E(x))_{x\in\Ind}$ of Banach spaces and $p\in[1,\infty]$, we can form the $\lp{p}$-direct sum of the $E(x)$ in the obvious way: this will be denoted by $\lpsum{p}_{i\in\Ind} E(x)$\/. 
Given a Banach algebra $A$\/, our definition of a Banach $A$-bimodule $M$ is the usual one: we require that the actions of $A$ on $M$ are jointly continuous, but not necessarily that they are contractive. When we write $M'$\/, we tacitly assume that it is equipped with the canonical $A$-bimodule structure obtained by taking adjoints of the actions of $A$ on $M$\/.

\subsection*{(Isometric) isomorphism of chain complexes and functors}
We assume familiarity with the notions of chain and cochain complexes of Banach spaces and modules.
For sake of brevity we adopt the convention that our chain and cochain complexes vanish in degrees $ \leq -1$\/, i.e.~are of the form
\[ 0 \larr E_0 \larr E_1 \larr \ldots \quad\text{ or }\quad
 0 \rarr M_0 \rarr M_1 \rarr \ldots \]

\begin{defn}
Let $A$ be a Banach algebra, and let $E_*$ and $F_*$ be chain complexes of left Banach $A$-modules. We say that $E_*$ and $F_*$ are \dt{topologically isomorphic as (module) chain complexes} if there exist mutually inverse chain maps $f: E_*\to F_*$ and $g:F_*\to E_*$\/, with each $f_n$ (and hence each $g_n$\/) a continuous $A$-module map.

If we can moreover arrange that each $f_n$ (and hence each $g_n$\/) is an isometry, we say that the chain complexes $E_*$ and $F_*$ are \dt{isometrically isomorphic}, and write $E_*\miso F_*$\/.
\end{defn}

\subsection*{Hochschild cohomology}
We repeat some background material in order to fix our notation.
Let $A$ be a Banach algebra and $M$ a Banach $A$-bimodule. The Hochschild cochain complex is
\begin{equation}\label{eq:cochaincomplex}
 0 \to \Co{C}{0}(A,M)\xrarr{\dif} \Co{C}{1}(A,M) \xrarr{\dif} \Co{C}{2}(A,M) \xrarr{\dif} \ldots
\end{equation}
where for each $n\in\Z_+$, $\Co{C}{n}(A,M)$ is the Banach space of all bounded $n$-linear maps from $A$ to $M$, and the coboundary operator $\dif:\Co{C}{n}(A,M)\to \Co{C}{n+1}(A,M)$ is given by
\[ (\dif \psi)(a_1,\dots a_{n+1}) \defeq
\left\{ \begin{aligned}
& a_1\psi(a_2,\dots, a_{n+1}) \\
+ & \sum_{j=1}^n (-1)^j
\psi(a_1,\dots, a_ja_{j+1}, \dots ,a_{n+1})
\\
+ & (-1)^{n+1} \psi(a_1,\dots, a_n)a_{n+1}
\end{aligned} \right.\]
(the proof that \eqref{eq:cochaincomplex} is a complex is straightforward).

We denote the kernel of $\dif: \Co{C}{n}(A,M)\to \Co{C}{n+1}(A,M)$ by $\Co{Z}{n}(A,M)$ and the range of $\dif:\Co{C}{n-1}(A,M)\to\Co{C}{n}(A,M)$ by $\Co{B}{n}(A,M)$. The quotient vector space $\Co{Z}{n}(A,M)/\Co{B}{n}(A,M)$ is the \dt{$n$th cohomology group of $A$} with \dt{coefficients in~$M$},
denoted by $\Co{H}{n}(A,M)$\/.

The case where $M=A'$ merits special attention. If $\Co{H}{n+1}(A,A')=0$ for all $n\geq 1$\/, we say that $A$ is \dt{simplicially trivial}.

\medskip
For most of this article $A$ will be the $\lp{1}$-convolution algebra of a discrete group~$G$\/.
There is a canonical one-dimensional $\lp{1}(G)$-module, denoted by $\Cplx_\veps$\/, corresponding to the augmentation character on $G$\/:
 we shall sometimes refer to $\Co{H}{n}(\lp{1}(G),\Cplx_\veps)$ as the \dt{$n$th bounded cohomology group of~$G$}.

Although we do not require much of the machinery of $\Ext$\, we shall assume familiarity with at least its basic definition and its relation to Hochschild cohomology, as can be found in~\cite[\S III.4]{Hel_HBTA}.
 Central to the machinery developed in \cite{Hel_HBTA} is the notion of an \dt{admissible} resolution or complex. We will need to consider a more precise notion.

\begin{defn}
Let $0\leftarrow E_0 \xlarr{d_0} E_1 \xlarr{d_1} \ldots$ be a chain complex
 of Banach spaces and continuous linear maps.
We say that the complex $E_*$ is \dt{$1$-split in $\Ban$} if there exist contractive linear maps $s_j:E_j\to E_{j+1}$\/, $j\geq 0$\/, such that $d_0s_0={\sf id}$ and
\[ s_{j-1}d_{j-1}+d_js_j ={\sf id}
 \qquad\text{ for all $j\geq 1$\/.} \]
\end{defn}

The point of introducing `$1$-splitness' explicitly is the following simple observation, whose proof we omit as it is straightforward.
\begin{lem}\label{l:sum-of-1split}
Let $\Ind$ be an index set and let $p\in[1,\infty]$. Suppose that for each $x\in\Ind$ we have a $1$-split chain complex
\[ 0\leftarrow E_0(x) \xlarr{d_0^x} E_1(x) \xlarr{d_1^x} \ldots \] 
in $\Ban$\/, such that for each $n$ we have $\sup_{x\in\Ind} \norm{d_n^x}<\infty$\/. Then the $\lp{p}$-sum
\[ 0\leftarrow \lpsum{p}_{x\in\Ind}E_0(x) \larr \lpsum{p}_{x\in\Ind}E_1(x) \larr \ldots \] 
is also a $1$-split complex, and is in particular exact.
\end{lem}

\begin{rem}
Without the $1$-split condition, the $\lp{p}$-sum of a family of exact chain complexes need not be exact: we shall return to this point in Section~\ref{s:statemain}.
\end{rem}
\end{section}

\begin{section}{Augmentation ideals}\label{s:statemain}
The original version of the following lemma was stated in the special case of augmentation ideals in discrete group algebras: the author thanks N.~\Gronbaek\ for pointing out that a more general result holds.

\begin{lem}\label{l:simptrivcondn}
Let $A$ be a unital Banach algebra which has a character $\varphi:A \to \Cplx$\/, and let $I=\ker(\varphi)$.
Then the following are equivalent:
\begin{itemize}
\item[$(i)$] $I$ is simplicially trivial;
\item[$(ii)$] $\Co{H}{n}(A,I')=0$ for all $n\geq 1$;
\item[$(iii)$] for each $n\geq 1$, the canonical map
\[\Co{H}{n}(A,\Cplx_\varphi) \xrarr{\varphi^*} \Co{H}{n}(A, A')\]
that is induced by the inclusion $\Cplx \to A', 1\mapsto \varphi$, is a topo\-logical isomorphism.
\end{itemize}
\end{lem}

\begin{proof}
The implications $(i) \iff (ii)$ are immediate from the observation that $A\iso\fu{I}$ and
the fact (see~\cite[Exercise III.4.10]{Hel_HBTA} or~\cite[Chapter 1]{BEJ_CIBA}) that $\Co{H}{n}(\fu{B},M)\iso\Co{H}{n}(B,M)$ for any Banach algebra $B$ and Banach $B$-bimodule $M$\/,
where $\fu{B}$ denotes the forced unitization of $B$\/.

To get the implications $(ii)\iff (iii)$, consider the long exact sequence of coho\-mology associated to the short exact sequence $0 \to \Cplx_{\varphi} \to A' \to I' \to 0$,~{\it viz}\/.
\[ \ldots \Co{H}{n}(A,\Cplx_{\varphi}) \xrarr{\varphi^*} \Co{H}{n}(A,A') \xrarr{\rho} \Co{H}{n}(A,I') \to \Co{H}{n+1}(A,\Cplx_\varphi) \to \ldots \]
We \emph{claim} that the map $\Co{H}{0}(A,A') \xrarr{\rho} \Co{H}{0}(A,I')$ is surjective. If this is true then our long exact sequence has the form
\[ 0\to \Co{H}{1}(A,\Cplx_{\varphi}) \xrarr{\varphi^*} \Co{H}{1}(A,A') \xrarr{\rho} \Co{H}{1}(A,I') \to \Co{H}{2}(A,\Cplx_\varphi) \to \ldots \]
and the equivalence of $(ii)$ and $(iii)$ now follows from \cite[Lemma 0.5.9]{Hel_HBTA}. Hence it remains only to justify our claim.

For any $A$-bimodule $X$, $\Co{H}{0}(A,X)$ is just the centre $Z(X)$ of $X$, so that $\rho:Z(A')\to Z(I')$ is given by the restriction of a trace on $A$ to the ideal~$I$. It therefore suffices to show that every element of $Z(I')$ extends to a trace on $A$. But this is easy: if $\psi \in I'$ and $\psi\cdot a =a\cdot\psi$ for all $a \in A$, then the functional $a\mapsto \psi(a-\varphi(a)\id[A])$ gives such a trace, and the proof is complete.
\end{proof}


We now specialize to group algebras. Throughout $G$ will denote a \emph{discrete} group, $\lp{1}(G)$ its convolution algebra and $\AI{G}$ the augmentation ideal in $\lp{1}(G)$ -- that is, the kernel of the augmentation character $\veps$ which sends each standard basis vector of $\lp{1}(G)$ to $1$.

\begin{defn}
A group $G$ is said to be \dt{\ct} if each element of $G\setminus\{\id[G]\}$ has an abelian centralizer.
\end{defn}
It is not immediately clear that there exist any non\-abelian, infinite, \ct\ groups: examples can be found in \cite[Chapter~1]{LynSchupp}, see in particular the remarks after Proposition~2.19. Let us just mention one family of examples.

\begin{thm}[Gromov, \cite{Grom_MSRI}; see also {\cite[Proposition 3.5]{Alo_WHGnotes}}]
Any torsion-free word-hyperbolic group is \ct.
\end{thm}
The arguments given for this in \cite{Grom_MSRI} are scattered over several sections and are not easily assembled into a proof. The simplest and clearest account appears in Chapter 3 of the survey article \cite{Alo_WHGnotes} (I would like to thank K.~Goda for drawing these notes to my attention).

\begin{rem}
It is often observed that direct products of hyperbolic groups need not be hyperbolic, the standard example being $F_2\times F_2$ where $F_2$ denotes the free group on two generators. In the current context it is worth pointing out that clearly $F_2\times F_2$ is not \ct\ (since the centralizer of $(\id, x)$ always contains a copy of $F_2\times\{\id\}$).
\end{rem}

\begin{thm}\label{t:augideal_coho_triv}
Let $G$ be a \ct, discrete group. Then for each $n \geq 1$, $\Co{H}{n}(\lp{1}(G), \AI{G}')=0$.
\end{thm}

The key to the proof is the following well-known idea: when we pass to a conjugation action, $\AI{G}$ decomposes as an $\lp{1}$-direct sum of modules of the form $\lp{1}(\Cl{x})$, where $\Cl{x}$ denotes the conjugacy class of $x$. Hence there is an isomorphism of cochain complexes
\begin{equation}\label{eq:decomp}
\Co{C}{*}(\lp{1}(G), \AI{G}')  \miso \lpsum{\infty}_{x \in \Ind} \Co{C}{*}(\lp{1}(G), \lp{1}(\Cl{x})' )
\end{equation}
where $\Ind$ is a set of representatives for each conjugacy class in $G\setminus\{\id[G]\}$. Our theorem will now follow from a computation of the coho\-mology of the complex on the right hand side of \eqref{eq:decomp}.

For each \emph{summand} on the right hand side of \eqref{eq:decomp}, the coho\-mology groups can be reduced to certain bounded coho\-mology groups: more precisely, it is observed 
in~\cite{PW_bddcoho} that for each $x$ there are isomorphisms
\begin{equation}\label{eq:pourabbas_white}
\begin{aligned}
\Co{H}{*}(\lp{1}(G), \lp{1}(\Cl{x})' )
 & \iso \Ext_{\lp{1}(G)}^*(\lp{1}(\Cl{x}),\Cplx) & \\
 & \iso \Ext_{\lp{1}(C_x)}^*(\Cplx,\Cplx) & \iso \Co{H}{*}(\lp{1}(C_x),\Cplx)
 \end{aligned}
\end{equation}
where $C_x$ denotes the centralizer of $x$. It is implicitly claimed in \cite[Corollary 3.7]{Pou_gpalg} that the coho\-mology of the cochain complex
\[ \lpsum{\infty}_{x \in \Ind} \Co{C}{*}(\lp{1}(G), \lp{1}(\Cl{x})' )\]
is isomorphic to
\[  \lpsum{\infty}_{x \in \Ind} \Co{H}{*}(\lp{1}(G), \lp{1}(\Cl{x})' ) \]
If this were the case then Theorem \ref{t:augideal_coho_triv} would follow immediately from Equation~\eqref{eq:pourabbas_white}. However,
the justification given in \cite{Pou_gpalg} for this supposed isomorphism is insufficient, because 
it is \emph{not} in general true that the coho\-mology of an {$\lp{\infty}$}-sum is the {$\lp{\infty}$}-sum of the coho\-mology of the summands.

As evidence we have the following simple example.

\begin{eg}\label{eg:SNIPER}
For each $n\in\Nat$ let $f_n:\Cplx\to\Cplx$ be the linear map `divide by $n$'\/: then the cokernel of each $f_n$ is zero, and so if we let $f: \lp{\infty}(\Nat)\to \lp{\infty}(\Nat)$ be the
$\lp{\infty}$-sum of all the~$f_n$ we find that
\[ \lpsum{\infty}_{n\in\Nat} \Coker f_n = 0 \neq \Coker f \;. \]
(It is precisely this phenomenon which motivates our somewhat laborious emphasis on $1$-split complexes and \emph{isometric} isomorphism of complexes.)
\end{eg}

\begin{rem}
In the special case where $G$ is \ct, each $C_x$ is abelian, hence amenable, and so for each $x$ the cochain complex $\Co{C}{*}(\lp{1}(C_x),\Cplx)$ has a contractive linear splitting.
Hence for such $G$\/, in order to deduce that the cochain complex
\[ \lpsum{\infty}_{x \in \Ind} \Co{C}{*}(\lp{1}(G), \lp{1}(C_x)' ) \]
splits, it would suffice to prove that the isomorphisms of Equation \eqref{eq:pourabbas_white} are induced by chain homotopies with norm control independent of $x$. This is implicitly done in \cite[\S4]{PW_bddcoho}, but only for second-degree coho\-mology.
\end{rem}

Rather than follow the approach outlined in the previous remark, we instead generalize
the argument sketched in the final section of~\cite{PW_bddcoho},
 so that it applies to any left $G$-set $S$ (i.e. we drop their hypothesis that the action is transitive). Since our hypotheses are weaker, we are not able to deduce isomorphism of coho\-mology groups as in \cite{PW_bddcoho}; however, our weaker conclusion suffices to prove Theorem \ref{t:augideal_coho_triv}.
\end{section}

\begin{section}{Disintegration over stabilizers}\label{s:mainhack}
The promised generalization goes as follows:
\begin{thm}\label{t:disintegrate_over_stabilizers}
Let $G$ be a discrete group acting from the left on a set $S$, and let $S=\coprod_{x \in \Ind} {\sf Orb}_x$ be the partition into $G$-orbits.
Regard $\lp{1}(S)$ as a Banach $\lp{1}(G)$-bimodule with left action given by the $G$-action on $S$ and right action given by the augmentation action $(x,g)\mapsto x$\/.

Let $H_x\defeq {\rm Stab}_G(x)$. Then for each $n$\/,
$\Co{H}{n}(\lp{1}(G), \lp{1}(S)')$ is
topologically isomorphic to the $n$th coho\-mology group of the complex
\[ 0 \rarr \lpsum{\infty}_{x \in\Ind} \Co{C}{0}(\lp{1}(H_x) ,\Cplx)
 \rarr \lpsum{\infty}_{x \in\Ind} \Co{C}{1}(\lp{1}(H_x) ,\Cplx)
 \rarr \ldots  \]
\end{thm}

\begin{coroll}\label{c:amenable_stabs}
Let $G$, $S$ be as above, and assume that each stabilizer subgroup $H_x$ is amenable. Then $\Co{H}{n}(\lp{1}(G), \lp{1}(S)')=0$ for all $n \geq 1$.
\end{coroll}
\begin{proof}[Proof of corollary]
Since each $H_x$ is amenable, the cochain complex $\Co{C}{*}(\lp{1}(H_x) ,\Cplx)$ admits a \emph{contractive} linear splitting in degrees $\geq 1$. Therefore the chain complex
\[ \lpsum{\infty}_{x \in\Ind} \Co{C}{*}(\lp{1}(H_x) ,\Cplx) \]
is also split in degrees $\geq 1$ by linear contractions, and is in particular exact in degree~$n$.
Now apply Theorem~\ref{t:disintegrate_over_stabilizers}.
\end{proof}


\begin{proof}[Proof of Theorem \ref{t:augideal_coho_triv}, assuming Corollary~\ref{c:amenable_stabs}]
By adapting the remarks preceding \cite[Theorem 2.5]{BEJ_CIBA}, it is straightforward to show that
\[ \Co{H}{n}(\lp{1}(G), \AI{G}') \iso \Co{H}{n}(\lp{1}(G),
\left(\twist{\AI{G}}\right)'\,) \]
where $\twist{\AI{G}}$ is the $\lp{1}(G)$-bimodule with underlying space $\AI{G}$ but with trivial right action and the conjugation left action.

(In more detail: there is a continuous chain isomorphism $\Theta^*$ from $\Co{C}{*}(\lp{1}(G),\AI{G}')$ to $\Co{C}{*}(\lp{1}(G),\left(\twist{\AI{G}}\right)')$, given by
\[ (\Theta^n\psi)(g_1,\ldots, g_n) = (g_1\dotsb g_n)^{-1}\cdot\psi(g_1,\ldots, g_n)\]
where $\psi\in\Co{C}{n}(\lp{1}(G),\AI{G}')$ and $g_1,\ldots, g_n\in G$\/. This formula differs slightly from those in \cite[\S2]{BEJ_CIBA}, because we wish to reduce to the case of cohomology coeff\-icients with augmentation action on the left, rather than on the right as in \cite{BEJ_CIBA}.
 The two viewpoints are essentially equivalent; but rather than convert between the two, it is simpler to verify that $\Theta^*$ is a chain map and that each $\Theta^n$ is an isomorphism of Banach spaces.)

Let $S= G \setminus \{ \id[G]\}$, regarded as a left $G$-set via conjugation action. Then the $\lp{1}(G)$-module $\twist{\lp{1}(G)}$ decomposes into a module-direct sum $\Cplx \oplus \lp{1}(S)$, where $\Cplx$ is the point module with trivial action. Composing the truncation map $\lp{1}(G) \to \lp{1}(S)$ with the inclusion map $\AI{G} \to \lp{1}(G)$ gives a linear isomorphism $\AI{G} \to \lp{1}(S)$, and this is also a $G$-module map (for the conjugation action). So for this action $\twist{\AI{G}} \cong \lp{1}(S)$ as $G$-modules, and therefore
\[ \Co{H}{n}(\lp{1}(G), \AI{G}') \iso \Co{H}{n}(\lp{1}(G),\lp{1}(S)') \]

Write $S$ as the disjoint union $S=\coprod_{x \in \Ind} \Cl{x}$ of conjugacy classes. The corresponding stabilizer subgroups are precisely the centralizers $C_x$ of each $x \in \Ind$; since $G$ is assumed to be \ct\ and $\id[G]\notin S$, each $C_x$ is commutative (hence amenable) and applying Corollary \ref{c:amenable_stabs} completes the proof.
\end{proof}

The proof of Theorem \ref{t:disintegrate_over_stabilizers} is broken into a succession of small lemmas: each is to some extent standard knowledge, but for our purposes we need to make explicit certain uniform bounds and linear splittings for which I can find no precise reference.
To do the requisite book-keeping, we take a functorial viewpoint.

\begin{notn}
The \dt{projective tensor product} of Banach spaces $E$ and $F$ will be denoted by $E\ptp F$\/.

Given a unital Banach algebra $B$\/, we denote by $\LunMod{B}$ the category whose objects are unit-linked, left Banach $B$-modules and whose morphisms are the $B$-module maps between them. $\Ban$ is the category of Banach spaces and bounded linear maps (equivalently, $\Ban\equiv \LunMod{\Cplx}$).

For such a $B$ there are two canonical functors: the `forgetful functor' $\cU:\LunMod{B}\to\Ban$\/, which sends a module to its underlying Banach space; and the `free functor' $B\ptp\blank :\Ban\to \LunMod{B}$\/, which sends a Banach space $E$ to the left $B$ module $B\ptp E$.

If $B$ is a Banach algebra, $M$ is a right Banach $B$-module and $N$ a left Banach $B$-module, we write $M\ptpR{B} N$ for the Banach tensor product of $M$ and $N$ over~$B$
(see \cite[\S II.3.1]{Hel_HBTA} for the definition and basic properties).
\end{notn}

\begin{defn}\label{dfn:iso-functor}
Let $B$\/, $C$ be unital Banach algebras, and let ${\mc F}$\/ and ${\mc G}$\/ be functors $\LunMod{B}\to\LunMod{C}$\/. We say that ${\mc F}$ and ${\mc G}$ are \dt{isometrically isomorphic} if there is a natural isomorphism $\alpha: {\mc F} \to {\mc G}$ such that, for each $M\in\LunMod{B}$\/, the morphism $\alpha_M:{\mc F}(M)\to {\mc G}(M)$ is an isometry as a map between Banach $C$-modules.
\end{defn}

\begin{rem}
Let $B$ be a Banach algebra and let $E_*$ be a chain complex in $\LunMod{B}$\/. If ${\mc F}$ and ${\mc G}$ are \emph{isometrically} isomorphic functors from $\LunMod{B}$ to $\Ban$\/, then the chain complexes ${\mc F}(E_*)$ and ${\mc G}(E_*)$ are isometrically isomorphic. In particular, if ${\mc F}(E_*)$ is $1$-split then so is~${\mc G}(E_*)$\/.
\end{rem}

\begin{lem}[Factorization of functors]\label{l:factorize_functors}
Let $B$ be a closed unital subalgebra of a unital Banach algebra $A$\/.
Regard $A$ as a right $B$-module via the inclusion homomorphism $B \hookrightarrow A$. Then:
\begin{itemize}
\item[$(i)$]
 we have a natural isometric isomorphism of functors
\[ A\ptpR{B}(B\ptp\blank) \miso A\ptp \blank \]
where $B\ptp\blank$ and $A\ptp \blank$ are the free functors from $\Ban$ (to $\LunMod{B}$ and $\LunMod{A}$ respectively;
\item[$(ii)$]
 we have a natural isometric isomorphism of functors
\[ \lHom{B}(\blank,\Cplx) \miso \lHom{A}\left( A\ptpR{B} \blank ,\,\Cplx\right)  \]
where both sides are functors  $\LunMod{B} \to \Ban$.
\end{itemize}
\end{lem}
The proof is clear (the analogous statements without the qualifier `isometric' are essentially given in, for instance, \cite[\S~II.5.3]{Hel_HBTA}).
\medskip

We shall also abuse notation slightly, to make some of the formulas more legible: if $H$ is a subgroup of $G$ and $M$ and $N$ are, respectively, right and left Banach $\lp{1}(H)$-modules, then we shall write $M\ptpR{H} N$ for the Banach tensor product of $M$ and $N$ over $\lp{1}(H)$\/.

\begin{lem}[A little more than flatness]\label{l:strongly_flat}
Let $H$ be any subgroup of $G$ and let $G/H=\{gH \st g \in G \}$ be the space of left cosets. Then we have a (natural) isometric iso\-mor\-phism of functors
\[ \cU_G \left(\lp{1}(G)\ptpR{H}\blank\right)  \miso \lp{1}(\Cst{G}{H}) \ptp ( \cU_H\blank) \]
where $\cU_G$ and $\cU_H$ are the forgetful functors
 to $\Ban$ (from the categories $\LunMod{\lp{1}(G)}$ and $\LunMod{\lp{1}(H)}$ respectively).
\end{lem}
\begin{proof}
Choose a transversal for $G/H$, that is, a function $\tau: G/H \to G$ such that $\tau({\mc J}) \in {\mc J}$ for all ${\mc J} \in G/H$. (Equivalently, $\tau({\mc J})H = {\mc J}$ for all ${\mc J}$). This transversal yields a function $\eta:G \to H$ such that
\[ g=\tau(gH)\cdot\eta(g)\quad\quad \text{ for all $g \in G$\/.} \]
Note that $\eta(gh)=\eta(g)\cdot h$ for every $g \in G$ and $h \in H$\/.

If $E$ is a unit-linked left $\lp{1}(H)$-module, define a contractive linear map $\lp{1}(G)\ptp E \to \lp{1}(G/H)\ptp E$ by $e_g \ptp v \mapsto e_{gH} \tp (\eta(g)\cdot v)$\/. This map factors through the quotient map $q: \lp{1}(G)\ptp E \to \lp{1}(G)\ptpR{H}E$, and so induces a linear contraction
\[ T_E: \lp{1}(G)\ptpR{H}E \to \lp{1}(\Cst{G}{H})\ptp E \]
where $T_E(e_g\tpR{H} v) \defeq e_{gH} \tp (\eta(g)\cdot v)$.

On the other hand, the composite map
\[ R_E: \lp{1}(\Cst{G}{H})\ptp E \xrarr{\tau\ptp{\sf id}_E}  \lp{1}(G)\ptp E \xrarr{q} \lp{1}(G)\ptpR{H}E  \]
is a linear contraction, defined by the formula $R(e_{\mc J}\tp v) \defeq e_{\tau({\mc J})} \tpR{H} v$. $R_E$ is the composition of two maps which are natural in $E$, hence is itself natural in $E$.
Direct checking on elementary tensors
shows that $R_E$ and $T_E$ are mutually inverse maps. Hence $R$ is a natural, isometric isomorphism from $\cU_G \left(\lp{1}(G)\ptpR{H}\blank\right)$ to $\lp{1}(G/H) \ptp ( \cU_H\blank)$ as required.
\end{proof}

\begin{lem}\label{l:reduce-to-one-sided-Ext}
Let $X$ be a left Banach $\lp{1}(G)$-module. Regard it as a $\lp{1}(G)$-bimodule $X_\veps$ by defining the right $G$-action on $X$ to be trivial (i.e.~augmen\-tation). Then for all $n$ there is a
topo\-logical iso\-morphism
\[ \Co{H}{n}(\lp{1}(G), X_\veps') \iso \Ext_{\lp{1}(G)}^n (X, \Cplx)\,. \]
\end{lem}
\begin{proof}
This is a special case of the isomorphisms
\[ \Co{H}{*}(A, \Lin{E}{F}) \iso \Ext_{A^e}^*(A,\Lin{E}{F}) \iso \Ext_A^*(E,F) \]
valid for any unital Banach algebra $A$ and any left Banach $A$-modules $E$ and $F$. (See \cite[Theorem III.4.12]{Hel_HBTA}.)
\end{proof}

\begin{lem}\label{l:prod_of_Hom}
Fix a Banach algebra $A$ and an index set $\Ind$\/; and for each $x\in\Ind$ let $0 \leftarrow M_0(x) \leftarrow M_1(x) \leftarrow \ldots$ be a chain complex of contractive left Banach $A$-modules and continuous $A$-module maps.

Suppose that for each $n\in\Nat$\/, the family of linear maps $(M_n(x)\to M_{n-1}(x))_{x\in\Ind}$ is uniformly bounded. Then for every left Banach $A$-module $N$\/, there is an isometric isomorphism of chain complexes
\[ \lHom{A}\left( \lpsum{1}_{x\in\Ind} M_*(x), N\right) \iso \lpsum{\infty}_{x\in\Ind} \lHom{A}(M_*(x), N) \]
\end{lem}
\begin{proof}[Outline of proof]

Let $n\geq 0$\/. Given $\psi : \lpsum{1}_x M_n(x) \to N$\/, define $\psi_y:M_n(y)\to N$ to be the map obtained by restricting $\psi$ to the embedded copy of~$M_n(y)$\/.
Then $(\psi_y)_{y\in\Ind}$ is a well-defined element of $\lpsum{\infty}_{y\in\Ind} \lHom{A}(M_n(y),
N)$\/. It is then straightforward to check that the function $\theta^n:\psi\mapsto (\psi_y)_{y\in\Ind}$ is an isometric linear isomorphism, and that the maps $\theta^n$ assemble to form a chain map.
\end{proof}

\begin{proof}[Proof of Theorem \ref{t:disintegrate_over_stabilizers}]
First observe that by Lemma \ref{l:reduce-to-one-sided-Ext}, there is a topological iso\-morphism
\[ \Co{H}{n}(\lp{1}(G), \lp{1}(S)') \iso \Ext^n_{\lp{1}(G)}(\lp{1}(S), \Cplx) \,. \]
Since $\Ext$ may be calculated up to topological isomorphism using any admissible projective resolution of the first variable, it therefore suffices to construct an admissible $\lp{1}(G)$-projective resolution $0 \leftarrow\lp{1}(S) \leftarrow P_0 \leftarrow P_1 \leftarrow \ldots$ with the following property:

\vspace{0.5em}
\noindent{${\mathbf (*)}$}\quad
\emph{the cochain complex}
$0 \rarr \lHom{\lp{1}(G)}(P_*,\Cplx)$
\emph{is topologically isomorphic to}
\[ 0 \rarr \lpsum{\infty}_{x \in\Ind} \Co{C}{*}(\lp{1}(H_x) ,\Cplx)  \;. \]
\vspace{0.5em}

We do this as follows. For each $x \in \Ind$, let
$0\leftarrow\Cplx\leftarrow P_*(x)$ denote the $1$-sided bar resolution of $\Cplx$ by left $\lp{1}(H_x)$-projective modules,~i.e.
\[  0\leftarrow  \Cplx \xlarr{\veps_x} \lp{1}(H_x) \xlarr{d^x_0} \lp{1}(H_x)^{\ptp 2} \xlarr{d^x_1} \ldots\]
where $\veps_x$ is the augmentation character and $d^x_n : \lp{1}(H_x)^{\ptp n+2} \to \lp{1}(H_x)^{\ptp n+1}$ is the $\lp{1}(H_x)$-module map given by
\[ d^x_n(e_{h(0)}\tp\ldots\tp e_{h(n+1)} ) = \left\{
\begin{gathered} \sum_{j=0}^n (-1)^j e_{h(0)}\tp\ldots\tp e_{h(j)h(j+1)}\tp\ldots\tp e_{h(n+1)} \\
+ (-1)^{n+1} e_{h(0)}\tp\ldots\tp e_{h(n)} \end{gathered}
\right. \]
for $h(0),h(1),\ldots, h(n+1)\in H_x$.
The complex $0\leftarrow\Cplx\leftarrow P_*(x)$ is $1$-split in $\Ban$
and so remains $1$-split after we apply $\lp{1}(G/H_x)\ptp\blank$ to it.
Therefore, by Lemma~\ref{l:strongly_flat} and the remark following Definition~\ref{dfn:iso-functor}, the chain complex
\begin{equation}\label{eq:induced_complex1}
0\leftarrow \lp{1}(\Cst{G}{H_x}) \xlarr{\tveps_x} \lp{1}(G)\ptpR{H_x}P_0(x) \xlarr{\td^x_0} \lp{1}(G)\ptpR{H_x}P_1(x) \xlarr{\td^x_1} \ldots
\end{equation}
is $1$-split as a complex in~$\Ban$\/.
Here, we have written $\tveps_x$ for the $\lp{1}(G)$-module map $\lp{1}(G)\tpR{H_x}\veps_x$\/, and $\td^x_n$ for the $\lp{1}(G)$-module map $\lp{1}(G)\tpR{H_x} d^x_n$\/.

For each~$n\geq 0$ let $P_n$ be the left Banach $\lp{1}(G)$-module
\[ P_n \defeq  \lpsum{1}_{x \in\Ind} \lp{1}(G)\ptpR{H_x} P_n(x) \quad; \]
write $\tveps$ for the $\lp{1}$-sum of all the $\tveps_x$\/, and define $\td_n$ similarly for each $n\geq 0$\/. 
As the $\lp{1}$-sum of $1$-split complexes is $1$-split (by Lemma~\ref{l:sum-of-1split}), the complex of Banach $\lp{1}(G)$-modules
\begin{equation}\label{eq:Good_Resoln}
0 \leftarrow \lpsum{1}_{x \in\Ind} \lp{1}(\Cst{G}{H_x})
 \xlarr{\tveps} P_0 \xlarr{\td_0} P_1 \xlarr{\td_1} \ldots
\end{equation}
is $1$-split as a complex in $\Ban$.
There is an isomorphism of $\lp{1}(G)$-modules
\[ \lp{1}(S)=\lp{1}\left(\coprod_{x \in \Ind} {\sf Orb}_x\right)\miso \lpsum{1}_{x \in \Ind} \lp{1}({\sf Orb}_x)\miso \lpsum{1}_{x \in\Ind} \lp{1}(\Cst{G}{H_x}) \]
where in the last step we identified the orbit of $x$ with the coset space $G/H_x$ via the correspondence $g\cdot x \leftrightarrow gH_x$.
Hence $0\leftarrow\lp{1}(S)\leftarrow P_*$ is an admissible complex of Banach $\lp{1}(G)$-modules.

Moreover, for each $x \in \Ind$ and $n \geq 0$,
 Lemma~\ref{l:factorize_functors} provides an isometric isomorphism of left $\lp{1}(G)$-modules
\[ \lp{1}(G)\ptpR{H_x} P_n(x) \miso \lp{1}(G)\ptp \lp{1}(H_x)^{\ptp n} \]
and taking the $\lp{1}$-direct sum over all $x$ yields isometric isomorphisms of left $\lp{1}(G)$-modules
\[ P_n =  \lpsum{1}_{x \in\Ind} \lp{1}(G)\ptpR{H_x} P_n(x)
 \miso  \lpsum{1}_{x \in\Ind}  \lp{1}(G)\ptp \lp{1}(H_x)^{\ptp n} \miso \lp{1}(G)\ptp \left(\lpsum{1}_{x \in\Ind}\lp{1}(H_x)^{\ptp n} \right)  \]
from which we see that each $P_n$ is free -- and hence projective -- as an $\lp{1}(G)$-module.

\emph{Combining the previous two paragraphs we see that $0\leftarrow\lp{1}(S)\leftarrow P_*$ is an admissible resolution of $\lp{1}(S)$ by $\lp{1}(G)$-projective modules.}

It remains to verify the condition {\bf ($*$)}. Observe that for each $x$
\[ \lHom{\lp{1}(H_x)}(P_*(x), \Cplx) \miso \Co{C}{*}(\lp{1}(H_x),\Cplx) \quad;\]
hence by Lemma \ref{l:factorize_functors} we have
\[  \Co{C}{*}(\lp{1}(H_x),\Cplx) \miso \lHom{\lp{1}(G)}\left(\lp{1}(G)\ptpR{H_x} P_*(x),\,\Cplx\right) \;,\]
and taking the $\lp{\infty}$-sum over all $x$ yields
\[  \begin{aligned}
\lpsum{\infty}_{x \in\Ind} \Co{C}{*}(\lp{1}(H_x),\Cplx)
  & \miso  \lpsum{\infty}_{x \in\Ind} \lHom{\lp{1}(G)}\left(\lp{1}(G)\ptpR{H_x} P_*(x),\,\Cplx\right) \\
 & \miso \lHom{\lp{1}(G)} ( P_*, \Cplx ) \end{aligned}  \]
where for the last isomorphism we appealed to Lemma \ref{l:prod_of_Hom}.
\end{proof}

\end{section}

\begin{section}{Corollaries and closing remarks}

\begin{coroll}\label{c:final_cor}
Let $G$ be a \ct, discrete group. Then $\AI{G}$ is simplicially trivial.
\end{coroll}
\begin{proof}
This is immediate from combining Lemma \ref{l:simptrivcondn} and Theorem \ref{t:augideal_coho_triv}.
\end{proof}

\begin{rem}
Recalling that biflat Banach algebras are simplicially trivial, it is natural to enquire if our result might follow from biflatness of $\AI{G}$.
To see that this is not always the case, observe that if $\AI{G}$ is biflat then $\Co{H}{2}(\AI{G},\Cplx_{\rm ann})=0$ by \cite[Theorem 4.13]{Sel_biflat}, while it is known that
\[ \Co{H}{2}(\AI{F_2},\Cplx_{\rm ann}) \iso \Co{H}{2}(\lp{1}(F_2),\Cplx) \neq 0 \;.\]
\end{rem}

While it is known that $\AI{G}$ is amenable if and only if $G$ is,
there appears to be no analogous characterization of precisely when $\AI{G}$ is biflat.

\begin{unsolved}
Let $G$ be a discrete group. If $I_0(G)$ is biflat, is $G$ amenable?
\end{unsolved}

\begin{rem}
We remarked earlier that $F_2\times F_2$ is not \ct.
The arguments above show that $\AI{F_2\times F_2}$ is \emph{not} simplicially trivial, since its second simplicial coho\-mology will contain a copy of the second bounded coho\-mology of $C_{(x,\id)}$ where $x \in F_2\setminus\{\id\}$. (To see that $\Co{H}{2}(\lp{1}(C_{(x,\id)}),\Cplx)$ is non-zero, observe that $C_{(x,\id)} \iso C_x\times F_2$ is the direct product of a commutative group with $F_2$, hence has the same bounded coho\-mology as $F_2$\,; by \cite[Proposition 2.8]{BEJ_CIBA} $\Co{H}{2}(\lp{1}(F_2),\Cplx)\neq 0$.)
\end{rem}

The question of what happens for augmentation ideals in \emph{non-discrete}, locally compact  groups is much trickier since measure-theoretic considerations come into play. Johnson and White have shown \cite{BEJ-MCW_AI} that the augmentation ideal of $PSL_2({\mathbb R})$
 is not even weakly amenable; in contrast, $PSL_2(\Z)$ is known to be \ct\ and so by our results its augmentation ideal is simplicially trivial.
\end{section}

\section*{Acknowledgements}
The results in this paper are taken, with some modifications in the presentation, from Chapter~2 of the author's PhD thesis \cite{YC_PhD}, which was done at the University of Newcastle upon Tyne with support from an EPSRC grant.
The chain complexes were drawn using Paul Taylor's {\tt diagrams.sty} macros.

Particular thanks are due to K.~Goda, N.~\Gronbaek, A.~Pourabbas and M.~C. White for useful exchanges.
The author also thanks the University of Manitoba for its support while this paper was written up.

\end{document}